\documentclass[12pt]{article}

\usepackage{amsmath,fullpage,graphicx,subcaption}
\usepackage{enumitem,overpic,xcolor}

\DeclareMathOperator{\Ei}{Ei}

\newcommand{\la}{\lambda}
\newcommand{\non}{\nonumber}
\newcommand{\eps}{\epsilon}
\newcommand{\beqa}{\begin{eqnarray}}
\newcommand{\eeqa}{\end{eqnarray}}
\newcommand{\beqas}{\begin{eqnarray*}}
\newcommand{\eeqas}{\end{eqnarray*}}
\newcommand{\beq}{\begin{equation}}
\newcommand{\eeq}{\end{equation}}

\newcommand{\ra}{\rightarrow}

\newcommand{\ee}{\mathrm{e}}
\newcommand{\ii}{\mathrm{i}}
\newcommand{\pdhfrac}[2]{\mathchoice{\frac{#1}{#2}}{#1/#2}{#1/#2}{#1/#2}}
\newcommand{\fdd}[2]{\pdhfrac{\mathrm{d}#1}{\mathrm{d}#2}}
\newcommand{\sdd}[2]{\pdhfrac{\mathrm{d}^2#1}{\mathrm{d}#2^2}}

\renewcommand{\d}{\mathrm{d}}

\begin{document}

\author{S. Jonathan Chapman\footnote{Mathematical Institute, AWB, ROQ,
    Woodstock Road Oxford OX2  6GG. jon.chapman@maths.ox.ac.uk}}

\title{Large-order perturbation theory of linear eigenvalue problems}

\maketitle
\begin{abstract}
We consider a class of linear eigenvalue problems depending on a small
parameter $\eps$ in which the series expansion for the eigenvalue  in
powers of $\eps$ is
divergent. We develop a new technique to determine the precise nature of
this divergence. We illustrate the technique through its  application to four
examples: the anharmonic oscillator, a simplified model of
equatorially-trapped Rossby
waves, and two simplified models based on quasinormal modes of
Reissner-Nordstr\"om-de Sitter black holes. 
\end{abstract}

\paragraph{Keywords} Exponential asymptotics; Stokes
phenomenon; divergent expansion; WKB theory; quantum mechanics.

\section{Introduction}

Linear eigenvalue problems of the form
\[
  L_\eps g = \la g,\]
where $L_\eps$ is a differential operator (depending on a parameter
$\eps>0$), are ubiquitous in applied mathematics and
theoretical physics. The eigenvalue $\la$ might correspond to an energy
level, the frequency of a normal mode of oscillation, or the growth
rate in a linear stability analysis.
Often $\eps$ is a small parameter, in which case it
is common to develop the perturbation series for the eigenvalue
$\la$ in powers of $\eps$,
\[ \la \sim \sum_{n=0}^\infty \eps^n \la_n.\]
    Sometimes this series diverges, and it is
    of interest to determine the nature of this divergence.
 From a purely numerical point of view understanding the large-order behaviour tells you how many terms to take in an optimal approximation, how large the smallest achievable error is, and  whether the series can be meaningfully resummed.
  However, as has long been known, large-order growth is controlled by non-perturbative effects that are not visible in ordinary perturbation theory.
  For eigenvalue problems, this can mean an exponentially small imaginary part of the eigenvalue, indicating an instability not visible in the regular perturbation expansion, or, in quantum mechanics, tunnelling between classically distinct regions,  or     exponentially-small level splittings.
       Here we provide a
straightforward approach to determining the precise asymptotic
behaviour of $\la_n$ for large $n$.  

One of the earliest examples of such a problem is the quantum
anharmonic oscillator \cite{Bender69,Bender73}, which we revisit in
Section \ref{ex3}. This problem has a long history; it
was the first
non-exactly-solvable problem tackled by the newly-written Schrodinger
equation in 1926, has practical applications ranging
from  quantum chemistry and atomic-molecular physics to crystal
lattice vibrations in solid-state theory, and serves as a simple
model for quantum field theory \cite{turbiner}. 
The seminal work by Bender and Wu \cite{Bender69,Bender73} established
nature of the divergence of the perturbation series for the ground
state energy. This work became the prototype for similar analyses in
many other quantum mechanical systems, in what is now known as
large-order perturbation theory \cite{Arteca90,Guillou90}.

The main technique of Bender and Wu is to analytically continue in
$\eps$ (typically until $\eps$ is negative), and then solve the
resulting problem by
combination of  Liouville-Green (WKB) and matched asymptotic approximations. This is a
delicate procedure, since the goal is to identify an exponentially
small component of $\la$ beyond-all-orders of the divergent asymptotic
series. Cauchy's integral formula is then used to determine the
coefficients in the  power series expansion of $\la$ on the positive real
$\eps$ axis  in terms of the values of
$\la$ on either side of the negative real $\eps$ axis. This technique
is ingenious, and has proved successful, but the details can be very
complicated. The present work aims to present an alternative
method, which we hope will be useful.

Intriguing also is the work of Dunne \& \"Unsal
\cite{DunneUnsal041701,DunneUnsal105009}, who use a uniform WKB approach for quantum mechanical systems with degenerate minima, focusing in particular on double-well
and sine-Gordon potentials.
They write the wave function in terms of a
scaled parabolic cylinder function
$\Psi(x) = D_\nu(u(x)/\sqrt{\eps})/\sqrt{u'(x)}$, where the index $\nu$ is a parameter to be determined along with $u(x)$. Although $u(x)$  and $\la$ can be determined locally in terms of $\nu$ by a regular perturbation expansion, $\nu$ itself is determined by imposing a global boundary condition (a symmetry condition at the midpoint between wells). This shows $\nu$ is exponentially close to an integer, and provides an exponential correction to $\la$. In \cite{DunneUnsal041701} this exponential correction is used to determine the late-order terms in the power series, in a similar way to Bender \& Wu.

In the applied mathematics literature an early example of such a
problem occurs in the work of Boyd and Natarov \cite{Boyd98}, who
consider a model problem for an equatorially-trapped Rossby wave in a
shear flow in the ocean or atmosphere. There the main interest is in
the imaginary part of the eigenvalue (corresponding to the growth rate
of instability)---the divergent perturbation series is purely real,
but there is an exponentially small imaginary part beyond all orders.
In \cite{Shelton} this problem is attacked in almost the reverse
direction to Bender and Wu---the divergent series is first found, and
used to determine the exponentially small imaginary component of the
eigenvalue via optimal truncation and Stokes phenomenon, rather than
the other way round. Simplifying and extending 
the procedure from 
\cite{Shelton} forms the basis of the present work.

We present our procedure through its application to four examples.
Each  follows the same general framework,
but the final part of the analysis differs slightly in each case. We hope this will allow the
interested reader to adapt the method to their own particular problem.

The general framework is:
\begin{description}
\item{(i) \textbf{Inner region}.} Each problem starts with a regular perturbation expansion. Typically the coefficients in the expansion of the eigenvalue are determined by imposing a regularity condition at the origin $x=0$, and each term in the expansion is a polynomial in $x$. This region determines the eigenvalue expansion completely, but it is difficult to extract the large-order behaviour from the resulting recurrence relations.
\item{(ii) \textbf{Outer region}.} The regular expansion in (i) is not uniform, and rearranges for large $x$, leading to a new expansion once $x = X/\eps$ has been rescaled. The outer problem is a singular perturbation problem, so that this expansion diverges in the usual form of factorial/power, driven by singularities away from $X=0$. The large-order behaviour of this divergence is easy to determine by now standard methods. In addition the late terms have an independent component driven by the divergent eigenvalue expansion.
  
\item{(iii) \textbf{Boundary layer in the late terms near $X = 0$}.}
  This is the new ingredient to the method. The key observation is that the large-order approximation of the outer expansion (ii) is also non-uniform, so that there is another inner region near $X=0$, now not in the small-$\eps$ expansion but in the large-order expansion. The resolution of this inner region links the two parts of the expansion in (ii), and determines the large-order behaviour of the eigenvalue.  
\end{description}


\section{Example 1: Simplified black holes}
\label{ex2}

We consider the model problem
\beq
2  (1-\eps x) (-\omega g+ xg') +    g    +(xg')'     =0, \qquad -\infty<x<0,\label{ex2maineqn}
\eeq
with
$g(0) = 1$ and $g(x) = o(\ee^{-x})$ as $x \ra -\infty$, where $0<\eps
\ll 1$. This is a much-simplified version of the problem in \cite{blackholes} concerning
quasinormal modes of  Reissner-Nordstr\"om-de Sitter black holes, keeping only those ingredients necessary to illustrate the methodology; very roughly speaking $g$ is the charged scalar field perturbation, $x$ is the radial distance measured from  the cosmological horizon, $\eps$ is the charge, and the
eigenvalue $\omega$ is the frequency of the mode.

As $x \ra -\infty$ the two possible behaviours are
\[ g(x) \sim \ee^{\eps x^2}, \qquad g(x) \sim x^{\omega},\]
while as $x \ra 0$ the two possible behaviours are
\[ g(x) \sim 1, \qquad g(x) \sim \log x.\]
The boundary conditions at $x=0$ and $x=-\infty$  each remove one
degree of freedom, so that there is a nonzero solution only if
$\omega$ takes particular values. The goal is to find the asymptotic
expansion of the eigenvalues,
\[ \omega \sim \sum_{n=0}^\infty \eps^n \omega_n,\]
as $\eps \ra 0$, and in particular the form of the divergence of
$\omega_n$ as $n \ra \infty$. To simplify the presentation we will focus on the leading eigenvalue, here and in each of our examples, but of course the methodology works for any eigenvalue.

\subsection{Inner region}
\label{sec2ex2}
We start with
 \[ 2  (1-\eps x) (-\omega g+ xg') +    g    +(xg')'  
   =0.
 \]
We expand
 \beq g = \sum_{n=0}^\infty  \eps^n g_n, \qquad
 \omega = \sum_{n=0}^\infty \eps^n \omega_n,\label{ex2exp}
 \eeq
to give at leading order
\[
 2  (-\omega_0 g_0 + x g_0') + g_0 + (x g_0')' = 0.
 \]
 The solution which is regular at the origin is
 \[ g_0 =  L_{\omega_0-1/2}(-2x),\]
 where $L_n(z)$ is the Laguerre function.
To avoid exponential growth as $x \ra -\infty$ we need the Laguerre function to be a polynomial, i.e.~we need $\omega_0-1/2$ to be a non-negative integer.
Choosing the first of these, $n=0$, gives the solution  $g_0 = 1$, $\omega_0 = 1/2$.
At next order
\[
 2   x g_1' + (x g_1')' - 2  \omega_1  =   -x .
\]
The  solution which is regular at $x=0$ and does not grow exponentially at minus infinity is
\[g_{1} = -\frac{x}{2}, \qquad \omega_1 = -\frac{1}{4}.\]
In general
\beq
 2   x g_n' + (x g_n')'- 2  \omega_n = -2  x (\omega_{n-1} g_0 +
\cdots + \omega_0 g_{n-1}  -x g_{n-1}') + 2  (\omega_{n-1}g_1 + \cdots
\omega_1 g_{n-1}),\label{gneqn} 
  \eeq
and the solution is of the form
\[ g_n = \sum_{i=1}^n a_{ni} x^i,\]
with
\beqa
 2  j a_{n,j} + (j+1)^2 a_{n,j+1} &=&- 2\sum_{k=j-1}^{n-1} \omega_{n-1-k}
 a_{k,j-1} + 2(j-1) a_{n-1,j-1} + 2  \sum_{k=j}^{n-1} \omega_{n-k}
 a_{k,j},\non \\ && \mbox{ } \hspace{4cm} \qquad\mbox{ for } j = 1,\ldots,n,\label{rec1} \\
 a_{n1} - 2 \omega_n & = & 0.\label{rec2}
\eeqa
We can iterate to find $\omega_n$ numerically. Figure \ref{fig1}(a)
shows $|\omega_n|^{1/n}$ as a function of $n$; the linear growth in $n$  is consistent with
factorial growth in $\omega_n$ at large $n$.
In principle we could extract the asymptotic behaviour as $n \ra
\infty$ from (\ref{rec1})-(\ref{rec2}), but this is not so
straightforward.
The method we now highlight determines $\omega_n$ for large $n$
without the need to analyse  (\ref{rec1})-(\ref{rec2}).

\subsection{Outer region}
\label{sec:2.2}
The expansion \eqref{ex2exp} is not uniform in $x$---it rearranges
when $x$ is large. In this section we develop the corresponding
expansion valid for large $x$.

To this end we set $\eps x = X$ to give
 \[ 2 (1-X) (-\omega g+ Xg') +  g    +\eps (Xg')'  
   =0.
 \]
Now expanding
\beq
g = \sum_{n=0}^\infty \eps^n g_n\label{gexp}
\eeq
gives, at leading order,
\[ 2  (1-X) (-\omega_0 g_0+ Xg_0') +  g_0      =0,
 \]
so that  
\[ g_0 = B(1-X)^{1/2}X^{\omega_0-1/2},\]
for some constant $B$.
For there to be no singularity at $X=0$ we require $\omega_0=1/2$, in
agreement with \S\ref{sec2ex2}.
To match with the inner expansion as $X \ra 0$ we require $g_0 \ra 1$ so that  $B=1$.
In general, equating coefficients of $\eps^n$, 
\beq
X( g_n + 2  (1-X)g_n')= - (X g_{n-1}')' + 2  (1-X) (\omega_1g_{n-1} + \cdots + \omega_n g_0).\label{ex2gneqn}
\eeq
We need to determine the late terms in the expansion, that is, the
behaviour of $g_n$ as $n \ra \infty$.
There are two sources of divergence in $g_n$: the usual
factorial/power divergence driven by differentiating $g_{n-1}$, and a
factorial/constant divergence driven by $\omega_n$.
For the first, we follow the usual procedure \cite{CKA} by
supposing that  
\beq
g_n \sim \frac{G\Gamma(n+\gamma)}{\chi^{n+\gamma}}
\eeq
as $n \ra \infty$, where $G$ and $\chi$ are functions of $x$ and
$\gamma$ is constant.
Then, equating coefficients of powers of $n$ gives, at leading order, 
\[  \chi'=2 (1-X). \]
Since this divergence is driven by the singularity in $g_0$ at $X=1$,
we have  $\chi(1) = 0$,
so that
\[ \chi = - (1-X)^2. \]
At next order we find
\[ (2-5X)G + 2X(1-X) G' =0,\]
giving
\[ G = \frac{\Lambda}{X(1-X)^{3/2}},\]
for some constant $\Lambda$.
Thus (absorbing $(-1)^{-\gamma}$ into $\Lambda$)  this part of $g_n$ satisfies
\[ g_n \sim \frac{\Lambda (-1)^n\Gamma(n+\gamma)}{X(1-X)^{3/2} (1-X)^{2n+2\gamma}}.\]
As $X \ra 1$,
\[ g_n \sim \frac{\Lambda (-1)^n \Gamma(n+\gamma)}{(1-X)^{3/2} (1-X)^{2n+2\gamma}}.\]
Comparing powers of $1-X$ with $g_0$ gives
\[ -\frac{3}{2}-2 \gamma = \frac{1}{2} \qquad \Rightarrow \qquad \gamma = -1,\]
so that 
\beq
g_n \sim \frac{\Lambda (-1)^n \Gamma(n-1)}{X (1-X)^{2n-1/2}}.\label{ex2out1}
\eeq
The other part of $g_n$, driven by the divergence of $\omega_n$,  is
given by $g_n \sim Q \omega_n$ where
\[X( Q + 2 (1-X) Q') = 2(1-X) g_0 = 2(1-X)^{3/2} ,\]
so that 
\[ Q =   (1-X)^{1/2}(\log X +C).\]
The presence of $\log X$ here means we need to modify slightly the ansatz $g_n \sim Q \omega_n$ (essentially we need $C$ to include a term proportional to $\log n$). If we set instead $g_n  = (Q_0 \log n + Q_1)\omega_n$  then
\beqas
X( Q_0 - 2 (1-X) Q_0') &=&0,\\
X( Q_1 - 2 (1-X) Q_1') &=&  2(1-X)^{3/2},
\eeqas
so that
\[ Q_0 \log n + Q_1 = (1-X)^{1/2}\left( C_0 \log n + C_1 +\log X\right).\]
Putting the two parts of $g_n$ together gives 
\beq
g_n \sim  \frac{\Lambda (-1)^n \Gamma(n-1)}{X(1-X)^{2n-1/2}}+
 (1-X)^{1/2}(C_0 \log n + C_1 +\log X)\omega_n.\label{ex2outer}
\eeq
To determine $\Lambda$ we need to match with an inner region in the
vicinity of the singularity at $X=1$.

\subsection{Inner region near $X=1$}
Motivated by both $g_0(X) = \sqrt{1-X}$ and by \eqref{ex2out1} we set $X = 1 - \eps^{1/2} \hat{x}$, $g = \eps^{1/4} \hat{g}$ to give
\beq
(1 -2   \hat{x}\eps^{1/2}\omega) \hat{g}+ (2   \hat{x}(-1+\eps^{1/2} \hat{x}) -\eps^{1/2})\hat{g}'  + (1-\eps^{1/2} \hat{x}) \hat{g}''  
=0.\label{ex2innereqn}
\eeq
In terms of  the inner variable
\beqa
g_0 & = & \eps^{1/4} \hat{x}^{1/2},\label{ex2g0in}\\
\eps^n g_n & \sim & \frac{\eps^{1/4}\Lambda (-1)^n \Gamma(n-1)}{\hat{x}^{2n-1/2}}.\label{ex2gnin}
\eeqa
At leading order in \eqref{ex2innereqn},
 \[\hat{g}_0 -  2\hat{x} \hat{g}_0' + \hat{g}_0'' = 0.\]
Writing
\[ \hat{g}_0 = \sum_{n=0}^{\infty}c_n  x^{1/2-2n},\]
we find
 \[ c_n = -\frac{(2n-5/2)(2n-3/2)c_{n-1}}{4n} , \qquad c_0=1,\]
 where the latter equation comes from matching with (\ref{ex2g0in}).
 Thus
\[ c_n = -(-1)^n\frac{(3/4)_{n-1} (5/4)_{n-1}}{16 (2)_{n-1}},\]
where $(x)_n = \Gamma(x+n)/\Gamma(x)$ is Pochammer's symbol.
Matching with \eqref{ex2gnin} gives
\[ \Lambda =\lim_{n \ra \infty} \frac{(-1)^n c_n}{\Gamma(n-1)} =-
  \frac{1}{16 \Gamma(3/4) \Gamma(5/4)} = -\frac{1}{4 \sqrt{2}\, \pi}.
\]

\subsection{Boundary layer in the late terms near $X = 0$}
\label{ex2:sec2.4} 
So far everything we have done has followed the standard approach to
finding the late terms of the expansion, as described in \cite{CKA}, for example.
In this section we make one crucial observation, which extends this
standard approach, and allows us link the two parts of the expansion
in (\ref{ex2outer}) and determine $\omega_n$.

This observation is that the large-$n$ asymptotic approximation for $g_n$ in
the outer region is non-uniform, and rearranges when $X$ is small.
We can see this directly from the asymptotic behaviour (\ref{ex2out1}), which is singular at
$X=0$, while we know that $g_n$ is in fact regular at $X=0$.

Thus there is another inner region near the origin, now not in the small-$\eps$ expansion of $g$, but in the large-$n$
expansion of $g_n$. To examine this inner region we rescale $X$ by
setting  $X = \xi/n$. Then the equation for $g_n$, equation \eqref{ex2gneqn}, 
becomes
\beq  \frac{\xi}{n} g_n + 2 \xi \left(1-\frac{\xi}{n}\right) g_n'=-n
  \left(\xi g_{n-1}'\right)' + 2  \left(1-\frac{\xi}{n}\right)
  (\omega_1g_{n-1} + \cdots + \omega_n g_0),\label{ex2finalin}
  \eeq
where $'$ is now $\d/\d\xi$.
Writing $X = \xi/n$ in \eqref{ex2outer}, the inner limit of the outer is
\beqa
g_n &\sim&  \frac{\Lambda (-1)^n n \Gamma(n-1)}{\xi(1-\xi/n)^{2n-1/2}}+
(1-\xi/n)^{1/2}(C_0 \log n + C_1 +\log \xi/n)\omega_n\non \\
& \sim &  \Lambda (-1)^n  \Gamma(n)\frac{\ee^{ 2 \xi}}{\xi}+
((C_0-1) \log n + C_1 +\log \xi)\omega_n.\label{ex2outin0}
\eeqa
This motivates writing
\[ \omega_n \sim \Omega (-1)^n \Gamma(n), \qquad
g_n \sim  H(\xi)\Omega (-1)^n \Gamma(n),\]
which, on substituting into \eqref{ex2finalin}, gives, at leading order,
\[ - (\xi H')'+  2 \xi H' =   2 .\] 
Thus
\beq
H =  \alpha_1 + \alpha_2 \Ei(2 \xi)+ \log \xi,\label{ex2inner0}
\eeq
where 
\[ \Ei(z) = \int_{-\infty}^z \frac{\ee^{t}}{t}\, \d t \]
is the exponential integral. Now $g_n$ should be regular as $\xi \ra 0$. Since
$\Ei(2 \xi) \sim \log \xi$ as $\xi \ra 0$,
we  need
\[\alpha_2 = -1\]
to remove the logarithmic singularity at $\xi=0$.
Since $\Ei(2\xi)$ exhibits Stokes phenomenon for large $\xi$ there
will be a switch in the  behaviour of the late terms depending on the
argument of $\xi$---this is what is known as the higher-order Stokes
phenomenon, a Stokes phenomenon not in the asymptotic expansion of $g$
as a function of $\eps$, but in the late-term approximation of
$g_n$ \cite{mortimer, howls}.
There is a higher-order Stokes line when $\xi$ crosses the
positive real axis, across which the constant contribution to the
large-$\xi$ approximation of $H$ (i.e. in the outer limit of the inner
expansion) changes. 
Note that there is no Stokes phenomenon associated with the particular
solution $\log \xi$, so that the coefficient of $\ee^{2 \xi}/\xi$ is fixed. This will
not be the case in our other examples.

To complete the analysis and determine $\Omega$ we need to match
\eqref{ex2inner0} with \eqref{ex2outin0}.
As $\xi \ra \infty$
\[ \Ei(2\xi) \sim \frac{\ee^{2\xi}}{2\xi}.\]
Matching with \eqref{ex2outin0} gives $C_0=1$ and $\Omega = -2\Lambda$
so that  
\beq
  \omega_n \sim -2 \Lambda (-1)^n \Gamma(n) = \frac{(-1)^{n}\Gamma(n)}{2
    \sqrt{2}\, \pi},\label{ex2omegan}
  \eeq
as $n \ra \infty$. 

\begin{figure}
  \begin{center}
\begin{subfigure}{0.4\textwidth}
   \begin{overpic}[width=\textwidth]{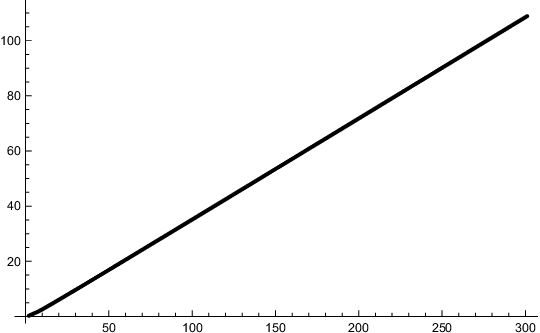}
      \put(-8,30){\rotatebox{90}{\scriptsize $\displaystyle |\omega_n|^{1/n}$}}
      \put(50,-4){\scriptsize $n$}
    \end{overpic}
    \vspace{-2mm}
    \caption{}
    \label{fig1a}
  \end{subfigure}\qquad \qquad
\begin{subfigure}{0.4\textwidth}
    \begin{overpic}[width=\textwidth]{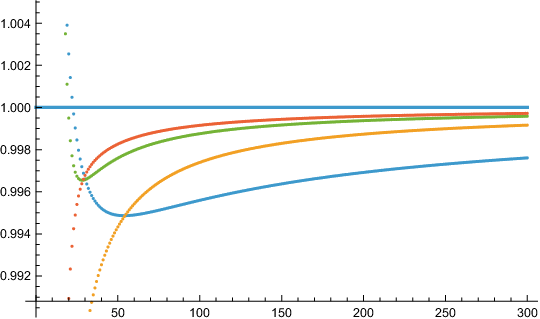}
      \put(-12,30){\rotatebox{90}{\scriptsize $\displaystyle2 \sqrt{2} \pi\frac{
          (-1)^n  \omega_n}{\Gamma(n)}$}}
      \put(50,-4){\scriptsize $n$}
    \end{overpic}
    \vspace{-2mm}
   \caption{}
  \label{fig1b}
\end{subfigure}
  \end{center}
  \caption{Divergence of the coefficients in the asymptotic expansion
    of $\omega$. (a) coefficients determined numerically  from
    (\ref{rec1})-(\ref{rec2}). The linear growth is consistent with
    factorial divergence. (b) The ratio of the numerical value to
  the asymptotic prediction (\ref{ex2omegan}). Blue is the base series, while
  orange,  green and red   correspond to enhanced convergence using
  Richardson extrapolation on two, three and four terms
  respectively. The black line is the asymptote, included to aid the
  eye. The convergence is slower than expected
  because of the presence of log terms in the higher-order corrections, unaccounted for in the extrapolation.}
  \label{fig1}
\end{figure}
In Fig.~\ref{fig1}(b) this result is compared with $\omega_n$ found by
numerically iterating (\ref{rec1})-(\ref{rec2}). The agreement is
found to be good, though the convergence is slightly slower than expected
  because of the presence of log terms in the higher-order corrections.

\section{Example 2: Anharmonic oscillator}
\label{ex3}
Having introduced the procedure with a simple model problem, we now
consider the classical problem of the anharmonic oscillator \cite{Bender69,Bender73}.
Much of the analysis follows the same framework, though the details of
the boundary layer in the late-term approximation analogous to \S
\ref{ex2:sec2.4} are a little different.

Consider\footnote{Note the typo in equation (1.1) of \cite{Bender73} in which the minus sign is missing.}
\[ \left(-\sdd{}{x} + \frac{x^2}{4} + \frac{\eps x^4}{4}\right) \Psi =
  \la \Psi,\]
with
\[ \Psi \ra 0 \qquad \mbox{ as } x \ra \pm \infty.\]
\subsection{Inner region}
We first factor out the decay at infinity by writing $\Psi = \ee^{-x^2/4} g$
to give
\beq
-g'' + x g' + \frac{g}{2} + \frac{\eps x^4 g}{4} - \la g = 0.
\eeq
Now expand
\beq
g = \sum_{n=0}^\infty \eps^n g_n, \qquad \la = \sum_{n=0}^\infty
\eps^n \la_n,\label{ex3exp}
\eeq
to give
\beqa
 -g_0'' + x g_0' + \frac{g_0}{2}  - \la_0 g_0 &=& 0,\\
 -g_n'' + x g_n' + \frac{g_n}{2}  - \la_0 g_n &=& - \frac{x^4 g_{n-1}}{4}
 + \sum_{k=1}^n \la_k g_{n-k}, \qquad n \geq 1.
\eeqa
For the first eigenvalue, the leading-order solution is $g_0=1$,
$\la_0=1/2$,  and, in general
\beq
g_n  =  \sum_{k=1}^{2n} a_{n,k} x^{2k}
\eeq
with
\beqa
2 k a_{n,k}& =& (2k+2)(2k+1) a_{n,k+1}-\frac{1}{4} a_{n-1,k-2} +
\sum_{i=1}^n \la_i a_{n-i,k},\qquad k = 2n, \ldots,1,\label{ex3rec1}\\
- 2 a_{n,1} & = & \la_n,\label{ex3rec2}
\eeqa
with the convention that $a_{n,k}=0$ for $k>2n$ and $k<1$. Equations
\eqref{ex3rec1}-\eqref{ex3rec2} are equivalent to eqn.~(6.3) in \cite{Bender73}.
It is argued in \cite{Bender73} that the leading-order late-term
behaviour of \eqref{ex3rec1}-\eqref{ex3rec2}  is the same as that of
the linearised equation (i.e. with the final sum omitted). With further
approximation, and quite a bit of analysis,  Bender \& Wu manage to extract the
leading-order behaviour of $\la_n$. Here we show how this may be
obtained by following the systematic procedure outlined in Section \ref{ex2}.

\subsection{Outer region}
As before, the expansion \eqref{ex3exp} is not uniform in $x$---it rearranges
when $x$ is large. In this section we develop the corresponding
expansion valid for large $x$.
We subtract off the leading-order eigenvalue by writing
\beq
\la = \frac{1}{2} + \eps \bar{\la}.
\eeq
We rescale into the far field by setting $\eps^{1/2} x = X$ to give
\beq
-\eps^2 g'' + \eps X g'+ \frac{X^4 g}{4 } - \eps^2 \bar{\la} g = 0.\label{ex3outeqn}
\eeq
The solution this time is of Liouville-Green (WKB) form,
\beq
g = \ee^{\phi/\eps}A, \qquad A \sim \sum_{n=0}^\infty \eps^n A_n.\label{ex3outexp}
\eeq
Substituting (\ref{ex3outexp}) into \eqref{ex3outeqn} and equating
coefficients  of
powers of $\eps$ gives, following a routine calculation,
\beq
\phi =
  \frac{1}{6}+ \frac{X^2}{4} - \frac{(1+X^2)^{3/2}}{6},
\qquad 
A_0 = \frac{\sqrt{2}}{(1+X^2)^{1/4} \sqrt{1+\sqrt{1+X^2}}},\label{ex3outA0}
\eeq
where the normalisation $\sqrt{2}$ comes from matching with the inner region.
In general, the equation for $A_n$ is
\begin{multline*}
  X(1+X^2)^{1/2} A_n' + \left((1+X^2)^{1/2} - \frac{1}{2} -
    \frac{1}{2(1+X^2)^{1/2}}\right)  A_n -  A_{n-1}''  \\
  -  \bar{\la}_0 A_{n-1} - \bar{\la}_1 A_{n-2} - \cdots -
  \bar{\la}_{n-1} A_0 = 0.
\end{multline*}
As before, there are two sources of divergence in $A_n$: the usual factorial/power
from repeated differentiation of the singularity in $A_0$, and a
factorial/constant divergence driven by  $\bar{\la}_n$.
For the first, we use the usual factorial/power ansatz   following the
procedure in \cite{CKA} to find
\beq
A_n \sim \frac{\Lambda}{(1+X^2)^{1/4} \sqrt{1-\sqrt{1+X^2}}
}\frac{ (-1)^n3^n \Gamma(n)}{(1+X^2)^{3n/2}}.\label{ex3outAn}
\eeq
The other part of $A_n$ satisfies
\[
  X(1+X^2)^{1/2} A_n' + \left((1+X^2)^{1/2} - \frac{1}{2} -
    \frac{1}{2(1+X^2)^{1/2}}\right)  A_n   \sim 
  \bar{\la}_{n-1} A_0,
\]
giving
\[ A_n \sim   \frac{\bar{\la}_{n-1}}{(1+X^2)^{1/4} \sqrt{1+\sqrt{1+X^2}}}
  \left(C_0\log n+C_1 -\tanh^{-1} \sqrt{1+X^2} \right),\]
where, as in \S\ref{sec:2.2}, the presence of a logarithm in $X$
necessitates the 
inclusion of a logarithm in $n$.
  Together
\begin{multline}
   A_n \sim
\frac{\Lambda}{(1+X^2)^{1/4} \sqrt{1-\sqrt{1+X^2}}
  }\frac{ (-1)^n3^n \Gamma(n)}{(1+X^2)^{3n/2}}\\+
  \frac{\bar{\la}_{n-1}}{(1+X^2)^{1/4} \sqrt{1+\sqrt{1+X^2}}}
  \left(C_0 \log n+C_1 -\tanh^{-1} \sqrt{1+X^2} \right).\label{ex3out}
\end{multline}
To determine $\Lambda$ we need to match with an inner region in the
vicinity of either $X=\ii$ or $X = -\ii$.
This problem is slightly unusual in that there are two singularities
in the leading-order solution, but they each produce a late-term
behaviour with the same singulant $\chi$, so that there is only one
factorial/power divergence in the late terms.

\subsection{Inner region near $X=\ii$}

We set $X = \ii - \ii\eps^{2/3} \hat{x}/2$, $A = \eps^{-1/6} \hat{A}$ to
give, at leading order,
\[ -2 \hat{x}^{1/2} \hat{A}' -
    \frac{ \hat{A}}{2 \hat{x}^{1/2}} +4 \hat{A}''  = 0.\] 
To match with \eqref{ex3outA0} requires
\beq  \hat{A} \sim \frac{\sqrt{2}}{ \hat{x}^{1/4}}
\qquad \mbox{ as }\hat{x} \ra \infty.\label{ex3A0inout}
\eeq
Writing
\beq
\hat{A} = \sqrt{2}\sum_{n=0}^{\infty}
\frac{g_n}{\hat{x}^{1/4+3n/2}},\label{ex3inout}
\eeq
gives
\[ g_n = -\frac{(6n-1)(6n-5)g_{n-1}}{12n} , \qquad g_0=1,\]
where the latter condition comes from \eqref{ex3A0inout}.
Thus
\beq
g_n = \frac{5 (-1)^n3^{n-2}(7/6)_{n-1} (11/6)_{n-1}}{4 (2)_{n-1}}.
\eeq
The inner limit of \eqref{ex3outAn} is
\beq
  \eps^n  A_n \sim
\frac{\Lambda}{\eps^{1/6}\hat{x}^{1/4}}\frac{ (-1)^n3^n
  \Gamma(n)}{\hat{x}^{3n/2}}.\label{ex3outin}
\eeq
Matching \eqref{ex3outin} with \eqref{ex3inout}  gives
\[ \Lambda =\sqrt{2}\lim_{n \ra \infty} \frac{g_n}{(-1)^n 3^n \Gamma(n)} =
  \frac{5\sqrt{2}}{36 \Gamma(7/6) \Gamma(11/6)} = \frac{1}{\pi \sqrt{2}}.
\]

\subsection{Boundary layer in the late terms near $X = 0$}

Here we come to the key step.   
The large $n$ asymptotic series for $A_n$ in the outer region
is nonuniform, and rearranges when $X$ is small---there is another inner region near the origin.
We emphasise again that this is a boundary layer not in the
small-$\eps$ expansion of $A$, but in the large $n$-expansion of
$A_n$. Again, this nonuniformity is evident because the asymptotic
formula \eqref{ex3out} is singular at $X=0$, while $A_n$ should be
regular  there.
This time the appropriate scaling for the inner region is $X =
\xi/n^{1/2}$ (so that $X A_n'$ balances with $A_{n-1}''$), giving
\beq
  \xi A_n' + \frac{3 \xi^2}{4n}  A_n -  n A_{n-1}''  +\cdots
  -  \bar{\la}_0 A_{n-1} - \bar{\la}_1 A_{n-2} - \cdots -
  \bar{\la}_{n-1} A_0 = 0.\label{ex3in}
\eeq
As $X \ra 0$ in \eqref{ex3out},
\beq
\eps^n  A_n \sim
\frac{\sqrt{2}\,\Lambda}{- \ii \sqrt{X^2}}
\frac{ (-1)^n3^n \Gamma(n)}{(1+X^2)^{3n/2}} +
\bar{\la}_{n-1}
\left(C_0\log n +C_1 \pm  \frac{\ii \pi}{2} - \log 2 + \log X
\right).
\label{ex3outin2a}
\eeq
Note that there are two choices of branch to be made here---one for
$\sqrt{X^2}$ arising from $\sqrt{1-\sqrt{1+X^2}}$  and one for the constant $\pm \ii \pi/2$ arising from
$\tanh^{-1}\sqrt{1+X^2}$.
In particular note that when matching to find $\Lambda$ we took $\sqrt{1-\sqrt{1+X^2}}$ to be real and
positive when $X$ approached $\pm \ii$, which means we need $-\ii
\sqrt{X^2}$ to be real and positive when $X$ is on the imaginary axis;
in turn this means we need $\sqrt{X^2} =X$ when $X$ is positive
imaginary, and  $\sqrt{X^2} =-X$ when $X$ is negative 
imaginary.
We will return to this choice and the
position of the branch cuts shortly, when we match with the inner
solution.

With $X = \xi/n^{1/2}$, \eqref{ex3outin2a} is
\beq
\eps^n  A_n \sim
\ii \sqrt{2}\,\Lambda n^{1/2}
(-1)^n3^n \Gamma(n)\frac{\ee^{-3\xi^2/2}}{\sqrt{\xi^2}}+
\bar{\la}_{n-1}
\left((C_0-1/2)\log n +C_1\pm \frac{\ii \pi}{2} - \log 2 + \log \xi
\right).
\label{ex3outin2}
\eeq
This motivates setting $A_n \sim H  \Omega(-1)^n 3^n \Gamma(n+1/2)$,
$\bar{\lambda}_{n-1} \sim \Omega (-1)^n 3^n \Gamma(n+1/2)$. Using this
ansatz in \eqref{ex3in} gives
\[   3\xi H' +  H''  - 3 =0,\]
so that
\beq
H =   \alpha_1 + \alpha_2 \int_{0}^\xi \ee^{-3u^2/2}\, \d u+3\int_0^\xi \ee^{-3t^2/2}\int_0^t
\ee^{3u^2/2}\, \d u \, \d t.\label{ex3H}
\eeq
Both the particular integral and the complementary function exhibit Stokes phenomenon for large $\xi$, so that there
will be a switch in the behaviour of the late terms depending on the
argument of $\xi$, corresponding to the higher-order Stokes phenomenon.
There is a higher-order Stokes line due to the particular integral when $\xi$ crosses the
real axis, across which the coefficient of $\ee^{-\xi^2}/\xi$ in
the far field (i.e. in the outer limit of the inner expansion)
changes. The branch cut associated with $\sqrt{\xi^2}$ in
\eqref{ex3outin2} must be chosen
to line up with this higher-order Stokes line.
In addition, there is a higher-order Stokes line on the imaginary axis,
across which the constant in the far field  changes. The branch cut
associated with $\pm \ii \pi/2$ must be chosen to align with this
higher-order Stokes line.
As $\xi \ra \infty$ in the first quadrant,
\beq
H \sim  \alpha_1 + \frac{\alpha_2 \sqrt{\pi}}{\sqrt{6}} +
\left(-\frac{\alpha_2}{3} - \frac{\ii \sqrt{\pi}}{\sqrt{6}} \right)
\frac{\ee^{-3\xi^2/2}}{\xi} + \cdots + \log \xi 
+ \frac{1}{2} (\gamma_E  +  \log 6) + \cdots.\label{ex3xiplus}
\eeq
where $\gamma_E$ is the Euler gamma.
As $\xi \ra \infty$ in the third quadrant,
\beq
H \sim \alpha_1 - \frac{\alpha_2 \sqrt{\pi}}{\sqrt{6}} +
\left(-\frac{\alpha_2}{3} + \frac{\ii \sqrt{\pi}}{\sqrt{6}} \right)
\frac{\ee^{-3\xi^2/2}}{\xi} + \cdots +\log \xi + 
  + \frac{1}{2} (\gamma_E  +  \log 6 - \ii \pi) + \cdots.\label{ex3ximinus}
\eeq
In the model problem of Section \ref{ex2} the key coefficient
$\alpha_2$ was determined by imposing that $H$ was regular at
$\xi=0$. In this case \eqref{ex3H} is regular at the origin for all
$\alpha_1$, $\alpha_2$, and it is matching with the outer solution
which determines $\alpha_2$. 

Matching \eqref{ex3xiplus} with \eqref{ex3outin2} as $\xi \ra \infty$
in the first quadrant gives
\[  \left(-\frac{\alpha_2}{3} - \frac{\ii \sqrt{\pi}}{\sqrt{6}} \right)  \Omega   =
 \ii \sqrt{2} \Lambda, \qquad C_0=1/2, \quad C_1 + \frac{\ii \pi}{2}
- \log 2 = \alpha_1 + \frac{\alpha_2 \sqrt{\pi}}{\sqrt{6}} +
\frac{\gamma_E}{2} + \log 6.\]
Matching \eqref{ex3ximinus} with \eqref{ex3outin2} as $\xi \ra \infty$
in the third quadrant  gives
\[  \left(-\frac{\alpha_2}{3} + \frac{\ii \sqrt{\pi}}{\sqrt{6}} \right)  \Omega  =
-  \ii \sqrt{2} \Lambda, \qquad C_0 =1/2, \quad
C_1 - \frac{\ii \pi}{2}
- \log 2 = \alpha_1 -\frac{\alpha_2 \sqrt{\pi}}{\sqrt{6}} +
\frac{\gamma_E}{2} + \log 6-\ii \pi
  .\]
Thus $\alpha_2=0$ and
\[ \Omega   =
-\frac{2 \sqrt{3}} { \sqrt{\pi}}  \Lambda  = -\frac{\sqrt{6}} { \pi^{3/2}} . \]
This gives, finally,
\beq
\bar{\lambda}_{n-1} = \la_n \sim  \frac{ (-1)^{n+1}\sqrt{6} } { \pi^{3/2}}\, 3^{n}\Gamma(n+1/2) ,\label{ex3lambdan}
\eeq
in agreement with \cite{Bender73}.
\begin{figure}
  \begin{center}
    \begin{overpic}[width=0.6\textwidth]{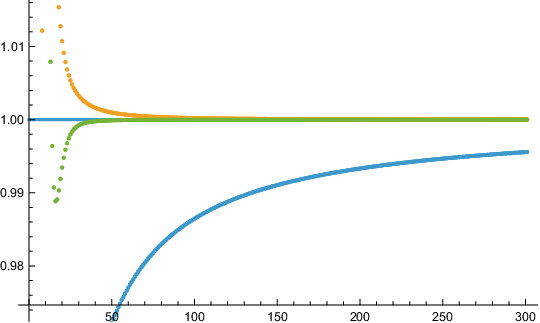}
      \put(-15,10){\rotatebox{90}{$\displaystyle \frac { \pi^{3/2} \la_n}{ (-1)^{n+1}\sqrt{6}\, 3^{n}\Gamma(n+1/2)  } $}}
      \put(50,-5){$n$}
    \end{overpic}
  \end{center}
 \caption{The ratio of the numerical value found by
 iterating (\ref{ex3rec1})-(\ref{ex3rec2}) to
  the asymptotic prediction (\ref{ex3lambdan}). Blue is the base series, while
  orange and green  correspond to enhanced convergence using
  Richardson extrapolation on two and three  terms respectively. The
  black line is the asymptote, included to aid the   eye. }
  \label{fig2}
\end{figure}
In Fig.~\ref{fig2} this result is compared with $\la_n$ found by
numerically iterating (\ref{ex3rec1})-(\ref{ex3rec2}); the agreement is
excellent.

\section{Example 3: Simplified Rossby waves}
\label{ex4}
Our third example is a simplified version of the model problem for an
equatorially-trapped Rossby wave considered in \cite{Shelton}.
Consider
\[
\sdd{\psi}{x} - 2 x \fdd{\psi}{x} + \frac{\eps^2 \psi}{1+\eps x} =
\la \psi,\qquad 
\ee^{-x^2/2}\psi \ra 0 \quad \mbox{ as } x \ra \pm \infty,
\]
with $\psi(0) = 1$. Essentially the problem considered in
\cite{Shelton} has the $\eps^2$ in the third term replaced with
$\eps$. The switch to $\eps^2$ makes the inner region below more complicated, but significantly
simplifies all the other regions of the analysis. This weakens the strength of the pole in Boyd and Natarov's Hermite-with-pole equation \cite{Boyd98}, while keeping the same structure.

\subsection{Inner region}
 \label{rossby:inner}
We will see that the expansion in powers of $\eps$ proceeds as 
\beq
  \psi = \sum_{n=0}^\infty \eps^n \psi_n, \qquad
  \la = \sum_{n=0}^\infty \eps^{2n} \la_n.\label{ex4:exp} 
\eeq
At leading order
\[ \psi_0'' - 2 x \psi_0' = \la_0 \psi_0.\]
In order for $\ee^{-x^2/2}\psi(y)$ to  decay as $x \ra \pm \infty$ we need $\psi_0$ to be a
Hermite polynomial.
The leading eigenvalue therefore has  $\psi_0 = 1$, $\la_0 = 0$.
In general
\[
  \psi_n'' - 2 x \psi_n' + \sum_{k=2}^n (-x)^{k-2}\psi_{n-k} =
  \sum_{k=0}^{\lfloor n/2\rfloor }\la_k \psi_{n-2k}, \qquad n \geq 1.\]
  In particular, we find $\psi_1 = 0$, while
 \[
  \psi_2'' - 2 x \psi_2' + 1  =
  \la_1 ,\] 
so that $\psi_2=0$, $\la_1 = 1$. Separating even and odd indices,
\beqas
\psi_{2n+1}'' - 2 x \psi_{2n+1}' + \sum_{k=1}^{n} (-x)^{2k-2}\psi_{2(n-k)+1} + \sum_{k=1}^{n} (-x)^{2k-1}\psi_{2(n-k)} &=&
\sum_{k=0}^{ n }\la_k \psi_{2(n-k)+1},\\
\psi_{2n}'' - 2 x \psi_{2n}' + \sum_{k=1}^{n} (-x)^{2k-2}\psi_{2(n-k)}  + \sum_{k=1}^{n-1} (-x)^{2k-1}\psi_{2(n-k)-1} &=&
\sum_{k=0}^{n }\la_k \psi_{2(n-k)}.
\eeqas
The solutions are of the form
\beq
\psi_{2n+1} = \sum_{k=1}^{n} a_{2n+1,k} x^{2k-1},\qquad
\psi_{2n} = \sum_{k=0}^{n-1} a_{2n,k} x^{2k},
\eeq
with
\beqa
2(2k-1)  a_{2n+1,k} & = & 2k(2k+1)a_{2n+1,k+1}
+ \sum_{m=1}^{n} a_{2(n-m)+1,k-m+1} \non\\ &&\qquad  \mbox{ }- \sum_{m=1}^{n} a_{2(n-m),k-m}
-
\sum_{m=0}^{ n }\la_m a_{2(n-m)+1,k} ,\label{ex4:rec1}\\
  4k  a_{2n,k}& = & (2k+2)(2k+1)a_{2n,k+1} +
\sum_{m=1}^{n} a_{2(n-m),k-m+1}  \non\\ && \qquad \mbox{ } - \sum_{m=1}^{n-1} a_{2(n-m)-1,k-m+1} -
\sum_{m=0}^{n }\la_m a_{2(n-m),k},\label{ex4:rec2}
\eeqa
with $a_{n,0}=0$ for $n>0$ and the convention that $a_{2n+1,k}=0$ and $a_{2n,k}=0$  if $k>n$ or $k<1$. For each $n$ equations \eqref{ex4:rec1}-\eqref{ex4:rec2} may be solved iteratively stepping down from $k=n$.
The solvability condition determining $\la_n$ comes from setting $k=0$ in \eqref{ex4:rec2}, giving $\la_n=2 a_{2n,1}$. We now follow the procedure of \S\ref{ex2} to determine the asymptotic behaviour of $\la_n$ as $n \ra \infty$.

\subsection{Outer region}
As usual, the expansion \eqref{ex4:exp} is not uniform in $x$ and rearranges when $x$ is large.
We set $\eps x = X$ to give
\[
\eps^2 \sdd{\psi}{X} - 2 X \fdd{\psi}{X} + \frac{\eps^2 \psi}{1+X} =
\la \psi.\]
The outer expansion now proceeds straightforwardly in powers of $\eps^2$ as
\beq
\psi = \sum_{n=0}^\infty \eps^{2n} \psi_n, \qquad \la = \sum_{n=0}^\infty \eps^{2n}\la_n.\label{ex4:outexp}
\eeq
At leading order
\[
- 2 X \fdd{\psi_0}{X} = \la_0 \psi_0,\]
with solution
\[ \psi_0 = B_0 X^{-\la_0/2}.\]
For there to be no singularity at $X=0$ requires $\la_0/2$ to be a non-positive integer, in agreement with the inner analysis in \S\ref{rossby:inner}. The leading eigenvalue therefore has $\la_0=0$, $\psi_0=1$.
At next order
\[
 - 2 X \fdd{\psi_1}{X} + \frac{1}{1+X} =
\la_1 ,\]
with solution
\[ \psi_1 = B_1 + \frac{(1-\la_1)}{2} \log X - \frac{1}{2} \log(1+X).\]
For there to be no singularity at $X=0$ we require $\la_1=1$, in agreement with \S\ref{rossby:inner}. The boundary condition $\psi_1(0) = 0$ (more properly a matching condition with the inner region) gives $B_1=0$.
In general
\beq
 \sdd{\psi_{n-1}}{X} - 2 X \fdd{\psi_n}{X} + \frac{ \psi_{n-1}}{1+X} =
 \sum_{k=1}^n \la_k \psi_{n-k}.\label{ex4:outn}
 \eeq
As usual, there are two types of divergence: a factorial/power from repeated  differentiation of the singularity $\log(1+X)$ in $\psi_1$, and a factorial/constant from  $\la_n$.
For the first, we use the usual ansatz to give
\beq
\psi_n \sim \frac{\Lambda \Gamma(n-1)}{X (1-X^2)^{n-1}}.\label{ex4part1}
\eeq
Here we see a curious feature of this example---the late term
behaviour $\psi_n$ was driven by a singularity at $X=-1$, but the
singulant vanishes also at $X=+1$.
Whereas in the anharmonic oscillator problem of \S\ref{ex3} both
singularities $X = \pm \ii$ were present in the early terms, here only $X=-1$ is
present in the early terms. The resolution of this apparent paradox,
as described in \cite{Shelton}, is a higher-order Stokes line which
turns off the contribution \eqref{ex4part1} in a region enclosing
$X=1$. We will return to this point later when
matching with an inner region near $X=0$.

The other part of $\psi_n$ satisfies $\psi_n \sim Q \la_n$ where
\[-2 X Q' = 1 ,\]
giving
\[ Q =   C-\frac{1}{2} \log X.\]
As usual, the presence of a logarithm  
means we need to  modify the large $n$ ansatz to essentially allow $C$ to depend on $\log n$. Putting both parts of $\psi_n$ together we have
\beq \psi_n \sim  \frac{\Lambda \Gamma(n-1)}{X (1-X^2)^{n-1}}+
\left(-\frac{1}{2} \log X + C_0 \log n + C_1\right)\la_n. \label{ex4:outer}
\eeq
The next step is to determine the constant $\Lambda$, by matching with an inner region near the singularity at $X=-1$.

\subsection{Inner region near $X=-1$}

We set $X= -1 + \eps^2 \hat{x}$.
Then the inner limit of the outer expansion satisfies
\beqa
\psi_0 + \eps^2 \psi_1 &\sim &  1 - \eps^2 \log \eps - \frac{\eps^2}{2} \log \hat{x} ,\label{ex4:outin0}\\
\eps^{2n}\psi_n& \sim& -\frac{\eps^2 \Lambda \Gamma(n-1)}{(2\hat{x})^{n-1}  }.\label{outin}
\eeqa
Equation (\ref{ex4:outin0}) motivates writing $\psi =  1 - \eps^2 \log \eps + \eps^2 \hat{\psi}$ to give the inner equation as
\[
\sdd{\hat{\psi}}{\hat{x}}
- 2 (-1+\eps^2 \hat{x}) \fdd{\hat{\psi}}{\hat{x}} + \frac{(1 - \eps^2 \log \eps + \eps^2 \hat{\psi})}{\hat{x}} = \la (1 - \eps^2 \log \eps + \eps^2 \hat{\psi}).
\]
At leading order
\[
\sdd{\hat{\psi}_0}{\hat{x}}
+ 2  \fdd{\hat{\psi}_0}{\hat{x}} + \frac{1}{\hat{x}} = 0.
\]
Thus
\[ \hat{\psi}_0 = \beta_1 + \beta_2 \ee^{- 2 \hat{x}} - \frac{1}{2} \log \hat{x} + \frac{1}{2} \ee^{- 2 \hat{x}} \Ei (2 \hat{x}) =  \beta_1 + \beta_2 \ee^{- 2 \hat{x}} - \frac{1}{2} \log \hat{x} + \frac{1}{2}\sum_{n=0}^\infty \frac{\Gamma(n+1)}{(2\hat{x})^{n+1}}.
\]
Matching with (\ref{outin}) gives
\[ \Lambda =-\frac{1}{2}.
\]

\subsection{Boundary layer in the late terms near $X = 0$}

As usual, the large $n$ asymptotic series for $\psi_n$ in the outer region
rearranges when $X$ is small.
This is again clear from the fact that the asymptotic approximation for $\psi_n$, \eqref{ex4:outer}, is singular at $X=0$, while $\psi_n$ is not.
There is a boundary layer near the origin  in the large $n$
expansion of $\psi_n$.
The appropriate scaling of this inner region is $X = \xi/n^{1/2}$, so that $\sdd{\psi_{n-1}}{X}$ balances $X \fdd{\psi_n}{X}$, so that \eqref{ex4:outn} becomes
\beq
n \sdd{\psi_{n-1}}{\xi} - 2 \xi \fdd{\psi_n}{\xi} + \frac{ \psi_{n-1}}{1+\xi/n^{1/2}} =
\sum_{k=1}^n \la_k \psi_{n-k}.\label{ex4:in2}
\eeq
The inner limit of the outer expansion \eqref{ex4:outer} is, for $X<0$,
\beqa
\psi_n &\sim& -\frac{\Gamma(n-1)n^{1/2}}{2\xi (1-\xi^2/n)^{n-1}}+
\left(-\frac{1}{2}\log (\xi) + C\right)\la_n \non \\
& \sim &-\frac{\Gamma(n-1/2)}{2}\frac{ \ee^{ \xi^2}}{\xi}+\left(-\frac{1}{2}\log (\xi) + C_0 \log n + C_1\right)\la_n.\label{sec4inout}
\eeqa
This motivates setting $\la_n = \Omega \Gamma(n-1/2)$, $\psi_n \sim  H \Omega \Gamma(n-1/2)$ in \eqref{ex4:in2}, to give
\[ H'' - 2 \xi H' = 1,\]
where $' \equiv \fdd{}{\xi}$, so that
\[ H = \alpha_1 + \alpha_2 \int_0^{\xi} \ee^{t^2}\, \d t + \int_0^{\xi} \ee^{t^2} \int_0^t \ee^{-p^2}\, \d p \, \d t.\]
Both the particular integral and the complementary function exhibit
Stokes phenomenon for large $\xi$, corresponding to the higher-order
Stokes phenomenon in $\psi$. 
There is a higher-order Stokes line due to the particular integral when $\xi$ crosses the
imaginary axis, across which the coefficient of $\ee^{\xi^2}/\xi$ in
the far field changes.
In addition, there is a higher-order Stokes line on the real axis,
across which the constant in the far field  changes. 
For real $\xi$, as $\xi \ra -\infty$
\[ H \sim \alpha_1 + \ii \frac{\alpha_2 \sqrt{\pi}}{2} +   \left( \frac{\alpha_2}{2} - \frac{\sqrt{\pi}}{4} \right) \frac{\ee^{\xi^2}}{\xi}+ \cdots - \frac{1}{4}\left(  \gamma_E + \log(-4\xi^2) + \cdots\right)
\]
while as $\xi \ra + \infty$
\[ H \sim \alpha_1 - \ii \frac{\alpha_2 \sqrt{\pi}}{2} +    \left( \frac{\alpha_2}{2} + \frac{\sqrt{\pi}}{4} \right) \frac{\ee^{\xi^2}}{\xi} + \cdots - \frac{1}{4}\left( \gamma_E + \log(-4\xi^2) + \cdots\right),
\]
where $\gamma_E$ is the Euler gamma.

\subsection{Matching with the outer}

\begin{figure}
  \begin{center}
    \begin{overpic}[width=0.6\textwidth]{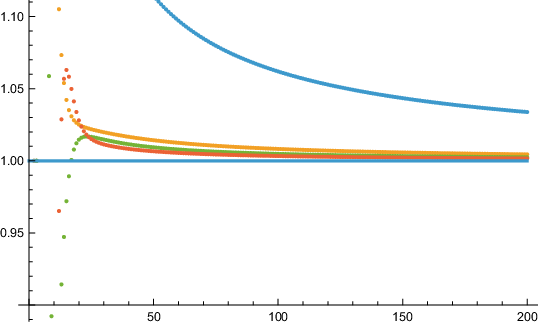}
      \put(-12,30){\rotatebox{90}{$\displaystyle \frac{\sqrt{\pi}
            \la_n}{\Gamma(n-1/2)}$}}
      \put(50,-5){$n$}
    \end{overpic}
  \end{center}
  \caption{The ratio of the numerical value found by iterating
    (\ref{ex4:rec1})-(\ref{ex4:rec2}) to
  the asymptotic prediction (\ref{ex4omegan}). Blue is the base series, while
  orange,  green and red   correspond to enhanced convergence using
  Richardson extrapolation on two, three and four terms respectively. The black line is the asymptote, included to aid the
  eye. The convergence is slower than expected
  because of the presence of logarithmic terms in the higher-order
  corrections, unaccounted for in the extrapolation.}
\label{fig3}
\end{figure}

The outer limit of the inner is
\beqas
 \psi_n &\sim&   \Omega \left( \frac{\alpha_2}{2} - \frac{\sqrt{\pi}}{4}\right) \frac{\ee^{\xi^2}}{\xi}\Gamma(n-1/2) + \cdots - \frac{\la_n}{2}\left( \log \xi  + \cdots\right)
\qquad \mbox{ as } \xi \ra -\infty,\\
\psi_n &\sim&   \Omega \left(  \frac{\alpha_2}{2} + \frac{\sqrt{\pi}}{4}\right) \frac{\ee^{\xi^2}}{\xi}\Gamma(n-1/2)+ \cdots - \frac{\la_n}{2}\left(
\log \xi + \cdots\right)
\qquad \mbox{ as } \xi \ra +\infty.
\eeqas
Matching with \eqref{sec4inout} as $\xi \ra -\infty$ gives
\[ \Omega \left(\frac{\alpha_2}{2} - \frac{\sqrt{\pi}}{4}\right) = - \frac{1}{2}.\]
For $X>0$, as per the discussion following \eqref{ex4part1}, there can
be no exponential term in the outer, because there must be no 
singularity at $X=1$. Thus, matching as $\xi \ra \infty$ gives
\[  \Omega\left(\frac{\alpha_2}{2} + \frac{\sqrt{\pi}}{4}\right) = 0.\]
Together
\[ \alpha_2 =  -\frac{\sqrt{\pi}}{2}, \qquad \Omega = \frac{1}{\sqrt{\pi}},\]
so that
\beq
\la_n \sim \frac{\Gamma(n-1/2)}{\sqrt{\pi}}.\label{ex4omegan}
\eeq
In Fig.~\ref{fig3} this result is compared with $\la_n$ found by
numerically iterating (\ref{ex4:rec1})-(\ref{ex4:rec2}); the agreement is
good, though the convergence is slower than expected
  because of the presence of logarithmic terms in the higher-order
  corrections.


\section{Example 4: Divergence driven by two singularities}

Our final example is chosen to illustrate that the divergence of the
eigenvalue can be driven by more than one singularity in the outer
solution, leading to more exotic behaviour. This is exactly what
happens in the model in \cite{blackholes} concerning
quasinormal modes of  Reissner-Nordstr\"om-de Sitter black holes.
The form of this divergence is more difficult to pick up with other
methods, and the interaction between two singularities makes it
difficult to guess the form of the divergence from numerical
calculations of the leading terms in the series.

Consider, as a model problem, 
\beq
\frac{b^2+(c+\eps x)^2}{b^2 + c(c+\eps x)} \left(-\omega g +  x
  g'\right)  + g +  (x g')'  = 0, \qquad -\infty<x<0, \label{ex5:eqn}
\eeq
with
     \[ g(0) = 1, \qquad g(x) = o(\ee^{-x}) \mbox{ as }x \ra \infty.\]
The relationship with \eqref{ex2maineqn} is clear---the coefficient of
the first term has been modified to generate an outer solution with
two singularities.

\subsection{Inner region}

\begin{figure}
  \begin{center}
    \begin{subfigure}{0.4\textwidth}
   \begin{overpic}[width=\textwidth]{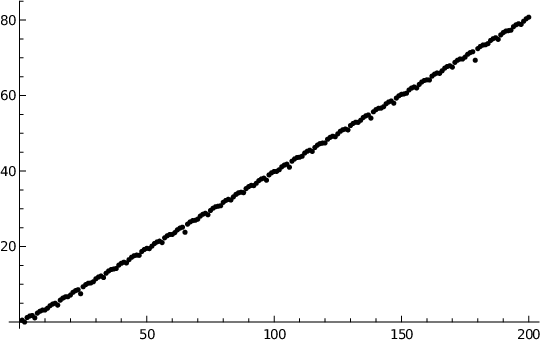}
      \put(-8,30){\rotatebox{90}{\scriptsize $\displaystyle |\omega_n|^{1/n}$}}
      \put(50,-4){\scriptsize $n$}
    \end{overpic}
    \vspace{-2mm}
    \caption{$b=1$, $c=1$}
    \label{fig4a}
  \end{subfigure}\qquad \qquad
\begin{subfigure}{0.4\textwidth}
    \begin{overpic}[width=\textwidth]{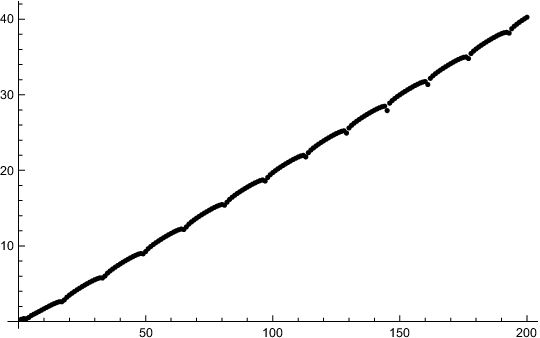}
      \put(-8,30){\rotatebox{90}{\scriptsize $\displaystyle |\omega_n|^{1/n}$}}
      \put(50,-4){\scriptsize $n$}
    \end{overpic}
    \vspace{-2mm}
    \caption{$b=1$, $c=3$}
  \label{fig4b}
\end{subfigure}
\end{center}
 \caption{Divergence of the coefficients in the asymptotic expansion
    of $\omega$, determined numerically  from
    (\ref{ex5:rec}). The linear growth is consistent with
    factorial divergence.}
  \label{fig4}
\end{figure}

We expand
 \[ g = \sum_{n=0}^\infty  \eps^n g_n, \qquad
   \omega = \sum_{n=0}^\infty \eps^n \omega_n,\]
to give at leading order
\[
  -\omega_0 g_0 + x g_0' + g_0 + (x g_0')' = 0,
\]
with solution  $g_0 = 1$, $\omega_0 = 1$.
At next order
\[ -\frac{c x}{b^2+c^2}  + \left(-\omega_1   +  x
    g_1'\right)   +  (x g_1')'  = 0,\]
with solution
\[g_{1} =  -\frac{c x}{b^2+c^2}, \qquad \omega_1 = \frac{c }{b^2+c^2}.\]
In general
\[ g_n = \sum_{i=1}^n a_{ni} x^i.\]
with
\begin{multline}
i a_{ni}=   \sum_{k=1}^n  \omega_k  a_{n-k,i}     -
 (i+1)^2  a_{n,i+1}  
  \\-  \frac{c}{ (b^2+c^2)}  \left(-2\sum_{k=0}^{n-1}  \omega_k  
  a_{n-k-1,i-1}    +2   (i-1) a_{n-1,i-1}  +
  a_{n-1,i-1}+  
 i^2  a_{n-1,i} \right) 
 \\ - \frac{1}{ (b^2+c^2)}\left(-\sum_{k=0}^{n-2}  \omega_k a_{n-2-k,i-2}
   +
  (i-2) a_{n-2,i-2}   \right),\label{ex5:rec}
\end{multline}
and $\omega_n = a_{n,1}$.
As usual, we can iterate \eqref{ex5:rec} numerically. Figure \ref{fig4}
shows $|\omega_n|^{1/n}$ as a function of $n$; the linear growth in $n$  is consistent with
factorial growth in $\omega_n$ at large $n$. Note that this growth is
not nearly as smooth as that in Fig.~\ref{fig1a}, with some ripples
present. Similar ripples can be seen in Fig.~2 of
\cite{blackholes}. These ripples are a direct result of the
interaction of the contributions from the two
singularities in the outer problem.

\subsection{Outer region}
We set $\eps x = X$ to give
 \[ \frac{b^2+(c+X)^2}{b^2 + c(c+X)} \left(-\omega g +  X
     g'\right)  + g +  \eps (X g')'  = 0.\]
Expanding
\beq
g = \sum_{n=0}^\infty \eps^n g_n, \qquad \omega = \sum_{n=0}^\infty
\eps^n \omega_n,\label{gexpEx4}
\eeq
and using $\omega_0=1$ gives
\[  \frac{b^2+(c+X)^2}{b^2 + c(c+X)} \left(- g_0 +  X
     g_0'\right)  + g_0       =0,
 \]
 with solution
 \beq
 g_0 = \frac{\sqrt{b^2+(c+X)^2}}{\sqrt{b^2+c^2}},\label{ex5:g0outer}
 \eeq
 where we have used the fact that $g_0 \ra 1$ as $X \ra 0$.
 We see that $g_0$ has singularities when $b^2+(c+X)^2=0$, i.e. $X = -c
 \pm \ii b$ (of course, the coefficient of the first term in \eqref{ex5:eqn}
 was chosen to make this the case).
In general
\[ \frac{b^2+(c+X)^2}{b^2 + c(c+X)}(- g_n + X g_n') + g_n =-(X
  g_{n-1}')' + \frac{b^2+(c+X)^2}{b^2 + c(c+X)} (\omega_1g_{n-1} + \cdots
  + \omega_n g_0).\] 
As usual there are two types of divergence: a factorial/power from the differentiation, and a factorial/constant from the $\omega_n$.
For the first, we use the usual ansatz 
\beq
g_n = \frac{G\Gamma(n+\gamma)}{\chi^{n+\gamma}}.
\eeq
At leading order in $n$ this gives
\[ -\frac{b^2+(c+X)^2}{b^2 + c(c+X)} = -   \chi'\]
so that
\[ \chi = \frac{(c+X)(c X + c^2-2b^2)}{2 c^2} + \frac{b^2(b^2+c^2)
    \log (b^2+c^2 + c X)}{c^3} + \mbox{const.} \]
This time there are two possible late term divergences, one
corresponding to each of the two singularities of the leading order
solution:
\beqas
\chi_1 & =& \frac{2 \ii b^3 + c(c+X)^2 - b^2(c+2X)}{2 c^2} + \frac{b^2(b^2+c^2)
 }{c^3} \log \left(\frac{b^2+c(c+X)}{b^2+ \ii b c}\right),\\
\chi_2 & =& \frac{-2 \ii b^3 + c(c+X)^2 - b^2(c+2X)}{2 c^2} + \frac{b^2(b^2+c^2)
}{c^3} \log \left(\frac{b^2+c(c+X)}{b^2- \ii b c}\right).
\eeqas
At next order,
\beqas
\frac{b^2+(c+X)^2}{b^2 + c(c+X)}(- G + X G') + G &=& G\chi' + X (2 G' \chi' + G \chi'').
\eeqas
Thus
\[ \frac{(b^4 + c(c+X)^2(c+3X) + b^2 (2 c^2 + 5 c X + 4 X^2))}{(b^2+c(c+X))^2} G + \frac{X (b^2 + (c+X)^2)}{b^2 + c(c+X)}  G' =0,\]
giving
\[ G = \frac{\Lambda (b^2+c(c+ X))}{X(b^2+ (c+X)^2)^{3/2}}.\]
Thus this part of $g_n$ satisfies
\[ g_n \sim  \frac{\Lambda_1 (b^2+c(c+ X))}{X(b^2+
    (c+X)^2)^{3/2}}\frac{\Gamma(n+\gamma_1)}{\chi_1^{n+\gamma_1}} + \frac{\Lambda_2 (b^2+c(c+ X))}{X(b^2+
    (c+X)^2)^{3/2}}\frac{\Gamma(n+\gamma_2)}{\chi_2^{n+\gamma_2}} .
\]
To determine $\gamma_1$ and $\gamma_2$ we match the order of the
singularity as $X \ra -c \pm \ii b$ with the early terms.
As $X \ra -c+ \ii b$,
\[ \chi_1 \sim \frac{(X + c - \ii b)^2}{c - \ii b},\]
\[ g_n \sim -\frac{\ii b}{2 \sqrt{2} (\ii b)^{3/2}(X + c - \ii
    b)^{3/2}}\frac{(c-\ii b)^{n+\gamma_1}}{(X + c - \ii b)^{2n+2\gamma_1}}
  \Lambda_1 \Gamma(n+\gamma_1).\]
Comparing powers of $X + c - \ii b$ with $g_0$ gives
\[ -\frac{3}{2}-2 \gamma_1 = \frac{1}{2} \qquad \Rightarrow \qquad
  \gamma_1 = -1.\]
A similar comparison as $X \ra -c-\ii b$ gives $\gamma_2=-1$ also, so that  
\[ g_n \sim  \frac{\Lambda_1 (b^2+c(c+ X))}{X(b^2+
    (c+X)^2)^{3/2}}\frac{\Gamma(n-1)}{\chi_1^{n-1}} + \frac{\Lambda_2 (b^2+c(c+ X))}{X(b^2+
    (c+X)^2)^{3/2}}\frac{\Gamma(n-1)}{\chi_2^{n-1}} .
\]
The other part of $g_n$ satisfies $g_n = (Q_0 \log n + Q_1)\omega_n$ where
\beqas
\frac{b^2+(c+X)^2}{b^2 + c(c+X)}(- Q_0 + X Q_0') + Q_0& =&
0,\\
\frac{b^2+(c+X)^2}{b^2 + c(c+X)}(- Q_1 + X Q_1') + Q_1& =&
\frac{b^2+(c+X)^2}{b^2 + c(c+X)}  g_0,
\eeqas
giving
\[ Q_0 \log n + Q_1 =   \frac{\sqrt{b^2+(c+X)^2}}{b^2 + c^2}\left(\log X
  +C_0 \log n + C_1\right).\]
Together 
\begin{multline}
  g_n \sim  \frac{\Lambda_1 (b^2+c(c+ X))}{X(b^2+
    (c+X)^2)^{3/2}}\frac{\Gamma(n-1)}{\chi_1^{n-1}} + \frac{\Lambda_2 (b^2+c(c+ X))}{X(b^2+
    (c+X)^2)^{3/2}}\frac{\Gamma(n-1)}{\chi_2^{n-1}}  \\+
\frac{\sqrt{b^2+(c+X)^2}}{b^2 + c^2}\left(\log X
  +C_0 \log n + C_1\right)\omega_n  .
\end{multline}
The next step is to determine $\Lambda_1$ and $\Lambda_2$ through
matching with inner regions near \mbox{$X = -c \pm \ii b$}.

\subsection{Inner region near $X=-c+ \ii b$}
To determine $\Lambda_1$ we look near $X=-c+ \ii b$.
We set $X = -c+ \ii b + \eps^{1/2}(c- \ii b)^{1/2} \hat{x}$, $g = \eps^{1/4} \hat{g}$ to give
\begin{multline*}
  \frac{ 2 \ii \eps^{1/2}(c- \ii b)^{1/2} \hat{x}}{b+ \ii c} \left(-\omega \hat{g}+ \frac{1}{(c- \ii b)^{1/2}\eps^{1/2}}(-c+ \ii b+\eps^{1/2}(c- \ii b)^{1/2} \hat{x})\hat{g}'\right) \\+  \hat{g}    +\eps \left(\frac{1}{(c- \ii b)^{1/2}\eps^{1/2}}\hat{g}' + (-c+ \ii b +(c- \ii b)^{1/2}\eps^{1/2} \hat{x}) \frac{1}{(c- \ii b)\eps}\hat{g}''\right)  
  =0.
\end{multline*}
At leading order
\[-  2\hat{x} \hat{g}_0' +\hat{g}_0 -\hat{g}_0'' = 0.\]
Writing
\beq \hat{g}_0 = \frac{\sqrt{2}(\ii b)^{1/2}
  (c-\ii
b)^{1/4}}{\sqrt{b^2+c^2}}\sum_{n=0}^{\infty}c_n  x^{1/2-2n},\label{ex5:outin}
\eeq
gives
\[ c_n = \frac{(2n-5/2)(2n-3/2)c_{n-1}}{4n} , \qquad c_0=1,\]
where the latter condition comes from matching with \eqref{ex5:g0outer}.
Thus
\[ c_n = -\frac{(3/4)_{n-1} (5/4)_{n-1}}{16 (2)_{n-1}}.\]
The inner limit of the outer expansion is
\beqa
\eps^n g_n &\sim& -\frac{\ii b \eps^n}{2 \sqrt{2} (\ii b)^{3/2}(X + c - \ii
  b)^{3/2}}\frac{(c-\ii b)^{n-1}}{(X + c - \ii b)^{2n-2}}
\Lambda_1 \Gamma(n-1)\non \\
&\sim& -\frac{\eps^{1/4}}{2 \sqrt{2} (\ii b)^{1/2}}\frac{\Lambda_1 \Gamma(n-1)\hat{x}^{1/2-2n}}{(c- \ii b)^{3/4}}.\label{ex5:inout}
\eeqa
Matching \eqref{ex5:outin} with \eqref{ex5:inout} gives
\beqas \Lambda_1 &=&- \frac{\sqrt{2}(\ii b)^{1/2} (c-\ii b)^{1/4}}{\sqrt{b^2+c^2}} 2 \sqrt{2} (\ii b)^{1/2}(c- \ii b)^{3/4} \lim_{n \ra \infty} \frac{c_n}{\Gamma(n-1)}\\
&=& \frac{(\ii b) (c-\ii b)^{1/2}}{(c+\ii b)^{1/2}} \frac{ 4}{16 \Gamma(3/4) \Gamma(5/4)} =\frac{(\ii b) (c-\ii b)^{1/2}}{\sqrt{2}\, \pi(c+\ii b)^{1/2}}.
\eeqas
A similar calculation near $X = - c - \ii b$ shows
\[ \Lambda_2 = \frac{(-\ii b) (c+\ii b)^{1/2}}{\sqrt{2}\, \pi(c-\ii b)^{1/2}}
 = \bar{\Lambda}_1,\]
where an overbar denotes complex conjugation.

\subsection{Boundary layer in the late terms near $X = 0$}
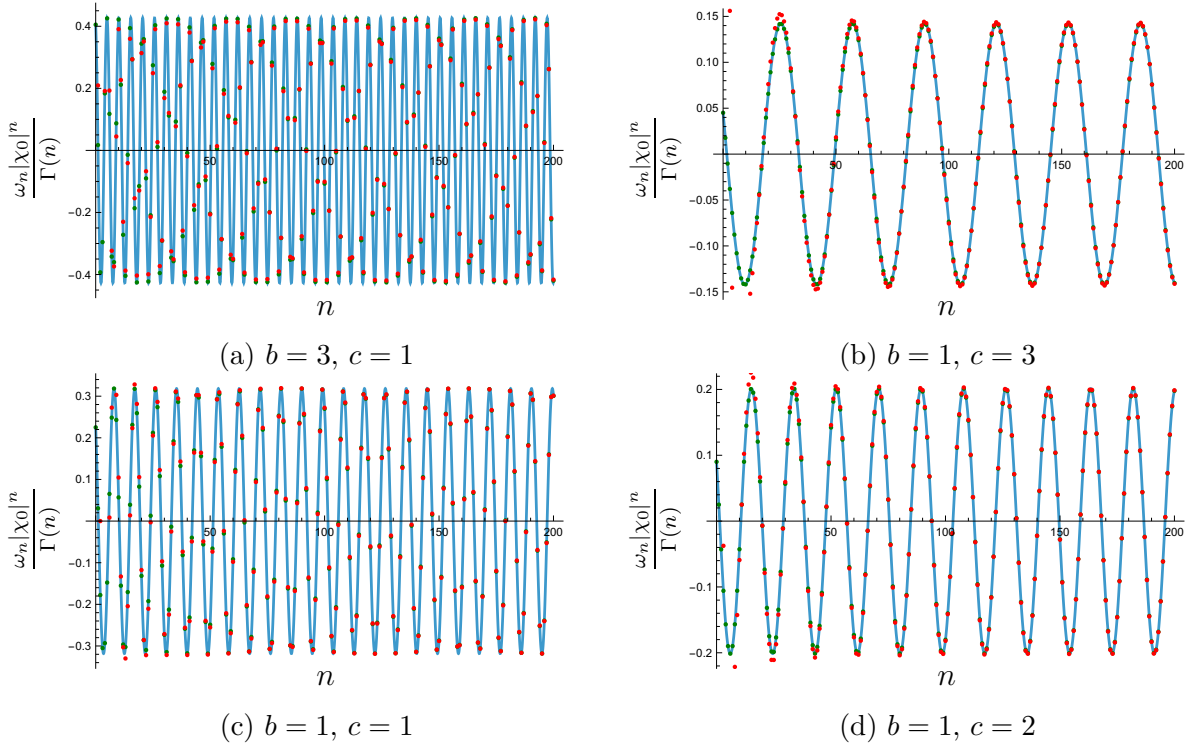
\begin{figure}
  \begin{center}
\begin{subfigure}{0.4\textwidth}
 \begin{overpic}[width=\textwidth]{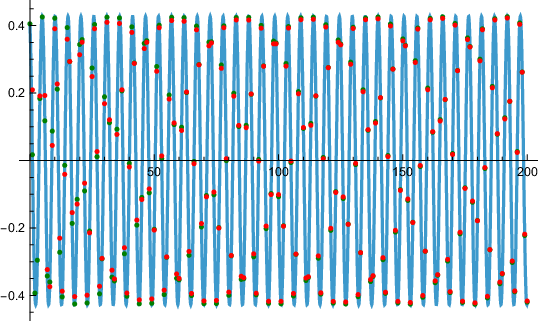}
      \put(-12,20){\scriptsize\rotatebox{90}{$\displaystyle
          \frac{\omega_n |\chi_0|^n}{\Gamma(n)}$}}
      \put(50,-3){$n$}
    \end{overpic}
    \vspace{-2mm}
  \caption{$b=3$, $c=1$}
    \label{fig5a}
 \end{subfigure}\qquad \qquad
\begin{subfigure}{0.4\textwidth}
 \begin{overpic}[width=\textwidth]{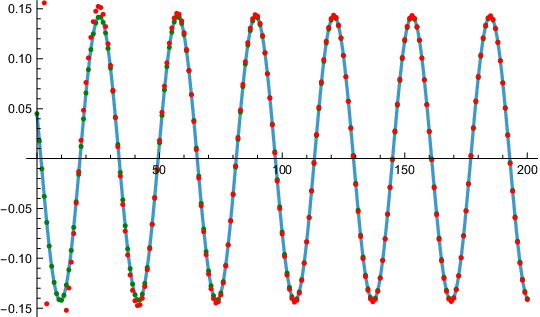}
      \put(-12,20){\scriptsize\rotatebox{90}{$\displaystyle
          \frac{\omega_n |\chi_0|^n}{\Gamma(n)}$}}
      \put(50,-3){$n$}
    \end{overpic}
    \vspace{-2mm}
  \caption{$b=1$, $c=3$}
    \label{fig5b}
 \end{subfigure}\\
\begin{subfigure}{0.4\textwidth}
 \begin{overpic}[width=\textwidth]{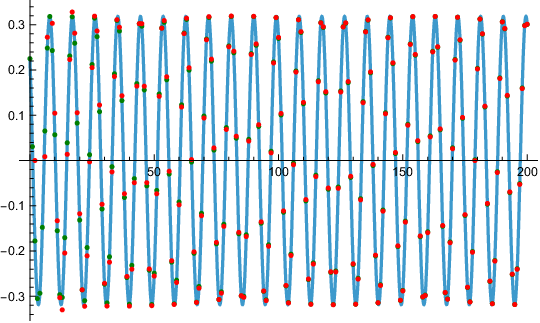}
      \put(-12,20){\scriptsize\rotatebox{90}{$\displaystyle
          \frac{\omega_n |\chi_0|^n}{\Gamma(n)}$}}
      \put(50,-3){$n$}
    \end{overpic}
    \vspace{-2mm}
  \caption{$b=1$, $c=1$}
    \label{fig5c}
\end{subfigure}\qquad \qquad
\begin{subfigure}{0.4\textwidth}
 \begin{overpic}[width=\textwidth]{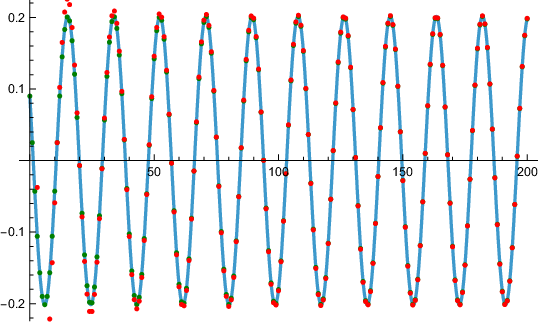}
      \put(-12,20){\scriptsize\rotatebox{90}{$\displaystyle
          \frac{\omega_n |\chi_0|^n}{\Gamma(n)}$}}
      \put(50,-3){$n$}
    \end{overpic}
     \vspace{-2mm}
   \caption{$b=1$, $c=2$}
    \label{fig5d}
 \end{subfigure}
  \end{center}
  \caption{A comparison of the asymptotic approximation
    \eqref{ex5:ans} with $\omega_n$ found by
numerically iterating \eqref{ex5:rec}, for various values of $b$ and
$c$. We normalise by $\Gamma(n)/|\chi_0|^{n}$ to remove the
exponential growth. The solid curve shows \eqref{ex5:ans} as a
continuous function of $n$, which is a sinusoidal oscillation of
period $\arg(\chi_0)/2\pi$. The green dots show \eqref{ex5:ans}
evaluated at integer $n$. The red dots are the numerical values.}
  \label{fig5}
\end{figure}

As in all our examples, the large $n$ asymptotic series for $g_n$ in the outer region
rearranges when $X$ is small, so that there is an inner region near the origin.
As in \S\ref{ex2} the appropriate rescaling is $X = \xi/n$, under
which the equation for $g_n$ 
becomes
\beq
\frac{b^2+(c+\xi/n)^2}{b^2 + c(c+\xi/n)}(- g_n + \xi g_n') + g_n
  =-n (\xi
  g_{n-1}')' + \frac{b^2+(c+\xi/n)^2}{b^2 + c(c+\xi/n)} (\omega_1g_{n-1} + \cdots
  + \omega_n g_0).\label{ex5:inner}
  \eeq
Writing $\chi_1$ and $\chi_2$ in terms of $\xi$ gives
\beqas
\chi_1 & \sim &  \frac{2 \ii b^3 - b^2 c + c^3}{2 c^2} +
\frac{b^2(b^2+c^2)}{c^3}\log\left(1 - \frac{\ii c}{b}\right) +
\frac{\xi}{n}+\cdots = \chi_0 + \frac{\xi}{n},\\
\chi_2 & \sim & \frac{-2 \ii b^3 - b^2 c + c^3}{2 c^2} +
\frac{b^2(b^2+c^2)}{c^3}\log\left(1 + \frac{\ii c}{b}\right) +
\frac{\xi}{n}+\cdots= \bar{\chi}_0 + \frac{\xi}{n},
\eeqas
say. Thus the inner limit of the outer solution is
\begin{eqnarray}
  g_n &\sim & \frac{\Lambda_1 }{\xi(b^2+
    c^2)^{1/2}}\frac{\Gamma(n)}{(\chi_0 + \xi/n)^{n-1}} + \frac{\bar{\Lambda}_1}{\xi(b^2+
    c^2)^{1/2}}\frac{\Gamma(n)}{(\bar{\chi_0} + \xi/n)^{n-1}}\non  \\
 && \hspace{6cm}  \mbox{ }+
\left(\log \xi/n
  +C_0 \log n + C_1\right)\omega_n\non  \\
&\sim & \frac{\Lambda_1 }{\xi(b^2+
    c^2)^{1/2}}\frac{\Gamma(n)}{\chi_0^{n-1}}\ee^{-\xi/\chi_0} + \frac{\bar{\Lambda}_1}{\xi(b^2+
    c^2)^{1/2}}\frac{\Gamma(n)}{\bar{\chi_0}^{n-1}}\ee^{-\xi/\bar{\chi}_0}\non  \\
 && \hspace{6cm} \mbox{ }+
\left(\log \xi
  +(C_0-1) \log n + C_1\right)\omega_n  \label{outin4}
\end{eqnarray}
From our analyses in \S\ref{ex2}-\S\ref{ex4} we have seen that the
boundary-layer approximation to $g_n$ comprises a particular integral 
driven by $\omega_n$ and a complementary function matching with the
remaining factorial/power divergence of the outer expansion.
For the current problem we write the particular integral as $g_n =H
\omega_n$, giving 
\[ (\xi H')'+  \xi H' =  1, \]
so that
\[ H = \log \xi\]
as in \S\ref{ex2}.
Matching this particular solution with \eqref{outin4} gives $C_0=1$,
$C_1 =  0$.  
The homogeneous solution may be written
\[ g_n = G(\xi)\frac{\Gamma(n)}{\chi_0^{n}}+
  \bar{G}(\xi)\frac{\Gamma(n)}{\bar{\chi}_0^{n}},\]
where

\[  \chi_0(\xi G')'+   \xi G'=0 ,\]
giving
\[ G = \alpha_1 + \alpha_2 \Ei(- \xi/\chi_0).\]
Thus together we have
\[ g_n \sim \left(\alpha_1 + \alpha_2 \Ei(- \xi/\chi_0)\right) \frac{\Gamma(n)}{\chi_0^{n}}
  +\left(\bar{\alpha}_1 + \bar{\alpha}_2 \Ei(-\xi/\bar{\chi}_0) \right) \frac{\Gamma(n)}{\bar{\chi}_0^{n}}
  +\omega_n \log \xi.\]
Now, as in \S\ref{ex2}, $g_n$ should be regular as $\xi \ra 0$. Thus,
since $ \Ei(\xi) \sim \log \xi$ as $\xi \ra 0$,
we  need
\[ \omega_n \sim  -\alpha_2 \frac{\Gamma(n)}{\chi_0^{n}} - \bar{\alpha}_2\frac{\Gamma(n)}{\bar{\chi}_0^{n}} .\]
To complete the analysis we need to match with \eqref{outin4} to determine $\alpha_2$.
As $\xi \ra \infty$
\[ \Ei(-\xi/\chi_0) \sim -\frac{\chi_0\ee^{-\xi/\chi_0}}{\xi}.\]
Thus the outer limit of the inner is
\[ g_n \sim \left(\alpha_1 - \alpha_2 \frac{\chi_0\ee^{-\xi/\chi_0}}{\xi}
    \right) \frac{\Gamma(n)}{\chi_0^{n-1}}
  +\left(\bar{\alpha}_1 - \bar{\alpha}_2\frac{\bar{\chi}_0\ee^{-\xi/\bar{\chi}_0}}{\xi}) \right) \frac{\Gamma(n)}{\bar{\chi}_0^{n-1}}
  +\omega_n \log \xi.\]
Matching with (\ref{outin4}) gives
\[ \alpha_2 = -\frac{\Lambda_1}{(b^2+c^2)^{1/2}}=
-\frac{\ii b}{\sqrt{2}\, \pi(c+\ii b)}
.\]
Thus
\beq
\omega_n \sim \frac{\ii b}{\sqrt{2}\, \pi(c+\ii b)} \frac{\Gamma(n)}{\chi_0^{n}}
  - \frac{\ii b}{\sqrt{2}\, \pi(c-\ii
    b)}\frac{\Gamma(n)}{\bar{\chi}_0^{n}} .\label{ex5:ans}
  \eeq
In Fig.~\ref{fig5} this result is compared with $\omega_n$ found by
numerically iterating \eqref{ex5:rec} for various values of $b$ and
$c$. The sinusoidal oscillation predicted by \eqref{ex5:ans} is clear
in Figs.~\ref{fig5b} and \ref{fig5d}, when the period of the oscillation
is long enough that there are many integers per cycle, but when the
period is short $\omega_n$ seems to jump around between different
longwave oscillations because of aliasing.

\section{Conclusion}

Through four examples, we have demonstrated a systematic procedure for
calculating the precise asymptotic behaviour of the late terms of the
asymptotic expansion of the eigenvalue in a variety of linear
eigenvalue problems. The framework in each of our examples is the
same.

After a regular perturbation expansion the eigenfunction at each order
is a polynomial, leading to a set of recurrence relations for the
coefficients of these polynomials and
the coefficients of the eigenvalue expansion.
While these relations are easy to iterate numerically to get the
leading terms of the eigenvalue expansion, it is hard to extract the
late term behaviour from them.

This regular perturbation expansion is nonuniform, and rearranges when
$x$ is large. Rescaling to an outer variable the corresponding outer
solution can be found, again as an asymptotic power series. This
series has a standard factorial/power divergence driven by
singularities in the leading-order approximation, and an additional
divergence driven by the divergent eigenvalue expansion. In contrast
to the original expansion, the late terms in this outer asymptotic
expansion are easy to find, using the usual factorial/power ansatz,
but the divergence of the eigenvalue is still undetermined.

However, the late term approximation of this outer expansion also
non-uniform, now not as $\eps \ra 0$ but as $n \ra \infty$. By
introducing a local variable in the equation for the late terms of
the outer expansion, a new inner expansion is generated in which the two
parts of the divergence become coupled, and the eigenvalue is
determined.

We hope our framework provides a template by which similar problems of
interest may be solved.

\subsubsection*{Acknowledgements}
I would like to acknowledge the Isaac Newton Institute for
Mathematical Sciences, Cambridge, for support and hospitality during the programme `Applicable resurgent asymptotics: towards a universal theory', where
part of the work on Examples 1 and 3 was undertaken, as well as useful
discussions with many of the participants. I would like to single out
in particular Ines Aniceto, who among other things  introduced me to reference
\cite{blackholes} during this programme.  This programme
was supported by EPSRC grant no. EP/R014604/1. 

I would also like to acknowledge the hospitality of the University of
Sydney and the support of Australian Research
Council Discovery Project  DP240101666.

For the purpose of Open Access, the author has applied a CC BY public copyright licence to any Author Accepted Manuscript (AAM) version arising from this submission.


\begin{thebibliography}{9}
\bibitem{Arteca90}
  Arteca, G. A.,  Fernández, F. M. \&   Castro, E. A., {\em Large Order
Perturbation Theory and Summation Methods in Quantum Mechanics},
Springer-Verlag, Berlin (1990).
\bibitem{Bender69} Bender, C.M, \& Wu, T.T., ``Anharmonic
  Oscillator,''
  {\em Phys. Rev.} {\bf 184}, 5, pp. 1231--1260 (1969).
\bibitem{Bender73} Bender, C.M, \& Wu, T.T., ``Anharmonic
  Oscillator. II. A Study of Perturbation   Theory in Large Order,''
  {\em Phys. Rev. D} {\bf 7}, 6, pp. 1620--1636 (1973).
\bibitem{Boyd98} Boyd, J. P. \& Natarov, A., ``A Sturm-Liouville eigenproblem of the fourth kind: A critical latitude
with equatorial trapping,'' {\em Stud. Appl. Math.} {\bf 101}(4),
433--455  (1998).
\bibitem{CKA} Chapman, S.J., King, J.R. \& Adams, K.L., ``Exponential
  Asymptotics and Stokes Lines in Nonlinear Ordinary Differential
  Equations,'' Proc. Roy. Soc. A {\bf 454}, pp. 2733--2755
  (1998).
\bibitem{mortimer} Chapman, S.J. \& Mortimer, D.B., ``Exponential
  asymptotics and Stokes lines in a partial differential equation,''
  Proc. Roy. Soc. A {\bf 461},  pp. 2385--2421   (2005). 
\bibitem{blackholes} Dias, O.J.C, Reall, H.S. \& Santos, J.E.,
  ``Strong cosmic censorship for charged de Sitter black holes with a
  charged scalar field,'' Class. Quantum Grav. {\bf 36}  045005 (2019).
\bibitem{DunneUnsal041701} Dunne, G.V. \& \"Unsal, M., ``Generating nonperturbative physics from perturbation theory,'' {\em Phys. Rev. D} {\bf 89}, 041701(R) (2014).
\bibitem{DunneUnsal105009} Dunne, G.V. \& \"Unsal, M., ``Uniform WKB, multi-instantons, and resurgent trans-series,'' {\em Phys. Rev. D} {\bf 89}, 1051009 (2014).
\bibitem{Guillou90}
 Le Guillou, J. C. \&  Zinn-Justin, J. (eds.), {\em Large-Order Behaviour of
Perturbation Theory} Elsevier Science Publishers, Amsterdam (1990).
\bibitem{howls} Howls, C.J., Langman, P.J. \& Olde Daalhuis, A.B.,
  ``On the higher order Stokes phenomenon,'' {\em Proc. R. Soc. A}
  {\bf 460}, pp. 1471–2946 (2005).
\bibitem{Shelton} Shelton, J., Chapman, S.J. \& Trinh, P.,
  ``Exponential Asymptotics for a Model Problem of an Equatorially
  Trapped Rossby Wave,'' {\em SIAM J. Appl. Math.} {\bf 84}, 4,
  pp. 1482--1503 (2024).
\bibitem{turbiner} Turbiner, A.V. \& del Valle Rosales, J.C., {\em Quantum Anharmonic Oscillator}, World Scientific, Singapore (2023).
\end{thebibliography}
\end{document}